\documentclass{amsart}
\usepackage{amsmath,amssymb,amsfonts}
\usepackage{eucal}

\usepackage{enumerate}

\theoremstyle{plain}
\newtheorem{theorem}{Theorem}[section]

\newtheorem{prop}[theorem]{Proposition}
\newtheorem{lemma}[theorem]{Lemma}

\theoremstyle{definition}
\newtheorem{definition}[theorem]{Definition}

\newtheorem*{thank}{Acknowledgments}

\input xy
\xyoption{all}

\newcommand{\Hom}{\text{Hom}}

\newcommand{\SSets}{\mathcal {SS}ets}
\newcommand{\sco}{\mathcal{SC}_\mathcal O}
\newcommand{\scxy}{\mathcal {SC}_{\{x,y\}}}

\begin{document}
\title[Simplicial Categories]%
{A Model Category Structure on the Category of Simplicial
Categories}

\author[J.E. Bergner]{Julia E. Bergner}

\address{Department of Mathematics \\
University of Notre Dame \\
Notre Dame, IN 46556}

\email{bergnerj@member.ams.org}

\date{\today}


\begin{abstract}
In this paper we put a cofibrantly generated model category
structure on the category of small simplicial categories.  The
weak equivalences are a simplicial analogue of the notion of
equivalence of categories.
\end{abstract}

\maketitle

\section{Introduction}
Simplicial categories, which in this paper we will take to mean
categories enriched over simplicial sets, arise in the study of
homotopy theories. Given any model category $\mathcal M$, the
simplicial localization of $\mathcal M$ as given in \cite{dk} is a
simplicial category which possesses the homotopy-theoretic
information contained in $\mathcal M$. Finding a model category
structure on the category of simplicial categories is then the
first step in studying the homotopy theory of homotopy theories.

In an early online version of the book \cite{dhk}, Dwyer,
Hirschhorn, and Kan present a cofibrantly generated model category
structure on the category of simplicial categories, but as
To\"{e}n and Vezzosi point out in their paper \cite{tv}, this
model category structure is incorrect, in that some of the
proposed generating acyclic cofibrations are not actually weak
equivalences.  Here we complete the work of \cite{dhk} by
describing a different set of generating acyclic cofibrations
which are in fact weak equivalences and which, along with the
generating cofibrations given in \cite{dhk}, enable us to prove
that the desired model category structure exists.

There are several contexts in which having such a model category
structure is helpful. In \cite{thesis}, we show that this model
structure is Quillen equivalent to three other model category
structures. These Quillen equivalences are of interest in homotopy
theory as well as higher category theory. To\"{e}n and Vezzosi
also use simplicial categories in their work on homotopical
algebraic geometry \cite{tv}.

Note that the term ``simplicial category" is potentially
confusing. As we have already stated, by a simplicial category we
mean a category enriched over simplicial sets.  If $a$ and $b$ are
objects in a simplicial category $\mathcal C$, then we denote by
$\Hom_\mathcal C (a,b)$ the function complex, or simplicial set of
maps $a \rightarrow b$ in $\mathcal C$.  This notion is more
restrictive than that of a simplicial object in the category of
categories. Using our definition, a simplicial category is
essentially a simplicial object in the category of categories
which satisfies the additional condition that all the simplicial
operators induce the identity map on the objects of the categories
involved \cite[2.1]{dk1}.

We will assume that our simplicial categories are small, namely,
that they have a set of objects.  A functor between two simplicial
categories $f:\mathcal C \rightarrow \mathcal D$ consists of a map
of sets $f:Ob(\mathcal C) \rightarrow Ob(\mathcal D)$ on the
objects of the two simplicial categories, and function complex
maps $f:\Hom_\mathcal C (a,b) \rightarrow \Hom_\mathcal D(fa,fb)$
which are compatible with composition.  Let $\mathcal {SC}$ denote
the category whose objects are the small simplicial categories and
whose morphisms are the functors between them. This category
$\mathcal {SC}$ is the underlying category of our model category
structure.

In a similar way, we can consider categories enriched over the
category of topological spaces. Making slight modifications to the
ideas from this paper, it is possible to put an analogous model
category structure on the category of small topological
categories.

Recall that a model category structure on a category $\mathcal C$
is a choice of three distinguished classes of morphisms, namely,
fibrations, cofibrations, and weak equivalences. We use the term
\emph{acyclic (co)fibration} to denote a map which is both a
(co)fibration and a weak equivalence.  This structure is required
to satisfy five axioms \cite[3.3]{ds}.

Before defining these three classes of morphisms in $\mathcal
{SC}$, we need some notation.  Suppose that $\mathcal C$ and
$\mathcal D$ are two simplicial categories.  Let $\pi_0 \mathcal
C$ denote the category of components of $\mathcal C$, namely, the
category in which the objects are the same as those of $\mathcal
C$ and the morphisms between objects $a$ and $b$ are given by
\[ \Hom_{\pi_0 \mathcal C}(a,b) = \pi_0
\Hom_\mathcal C(a,b). \] If $f:\mathcal C \rightarrow \mathcal D$
is a map of simplicial categories, then $\pi_0f: \pi_0 \mathcal C
\rightarrow \pi_0 \mathcal D$ denotes the induced map on the
categories of components of $\mathcal C$ and $\mathcal D$.

If $\mathcal C$ is a simplicial category, say that a morphism $e
\in \Hom_\mathcal C(a,b)_0$ is a \emph{homotopy equivalence} if it
becomes an isomorphism in $\pi_0 \mathcal C$.

Now, given these definitions, our three classes of morphisms are
defined as follows.

\begin{enumerate}
\item The weak equivalences are the maps $f:\mathcal C \rightarrow
\mathcal D$ satisfying the following two conditions:
\begin{itemize}
\item (W1) For any objects $a_1$ and $a_2$ in $\mathcal C$, the
map
\[ \Hom_\mathcal C (a_1,a_2) \rightarrow \Hom_\mathcal D (fa_1,fa_2) \]
is a weak equivalence of simplicial sets.

\item (W2) The induced functor $\pi_0f:\pi_0 \mathcal C
\rightarrow \pi_0 \mathcal D$ is an equivalence of categories.
\end{itemize}

\item The fibrations are the maps $f: \mathcal C \rightarrow
\mathcal D$ satisfying the following two conditions:
\begin{itemize}
\item (F1) For any objects $a_1$ and $a_2$ in $\mathcal C$, the
map
\[ \Hom_\mathcal C (a_1,a_2) \rightarrow \Hom_\mathcal D (fa_1,fa_2) \]
is a fibration of simplicial sets.

\item (F2) For any object $a_1$ in $\mathcal C$, $b$ in $\mathcal
D$, and homotopy equivalence $e:fa_1 \rightarrow b$ in $\mathcal
D$, there is an object $a_2$ in $\mathcal C$ and homotopy
equivalence $d:a_1 \rightarrow a_2$ in $\mathcal C$ such that
$fd=e$.
\end{itemize}

\item The cofibrations are the maps which have the left lifting
property with respect to the maps which are both fibrations and
weak equivalences.
\end{enumerate}

The weak equivalences are often called \emph{DK-equivalences}, as
they were first described by Dwyer and Kan in \cite{dk1}.  They
are a generalization of the notion of equivalence of categories to
the simplicial setting.

We now state our main theorem.

\begin{theorem} \label{MC}
There is a model category structure on the category $\mathcal
{SC}$ of all small simplicial categories with the above weak
equivalences, fibrations, and cofibrations.
\end{theorem}

We will actually prove the stronger statement that the above model
category structure is cofibrantly generated. Recall that a
cofibrantly generated model category is one for which there are
two specified sets of morphisms, one of generating cofibrations
and one of generating acyclic cofibrations, such that a map is a
fibration if and only if it has the right lifting property with
respect to the generating acyclic cofibrations, and a map is an
acyclic fibration if and only if it has the right lifting property
with respect to the generating cofibrations.  For more details
about cofibrantly generated model category structures, see
\cite[Ch. 11]{hirsch}. To prove the theorem, we will use the
following proposition, which is stated in more general form by
Hirschhorn \cite[11.3.1]{hirsch}.

\begin{prop} \label{CofGen}
Let $\mathcal M$ be a category with all small limits and colimits
and specified classes of weak equivalences and fibrations. Define
a map to be a cofibration if it has the left lifting property with
respect to the acyclic fibrations. Suppose that the class of weak
equivalences is closed under retracts and satisfies the ``two out
of three property" \cite[3.3]{ds}. Suppose further that there
exist sets $C$ and $A$ of maps in $\mathcal M$ satisfying the
following properties:
\begin{enumerate}
\item Both $C$ and $A$ permit the small object argument
\cite[10.5.15]{hirsch}.

\item A map is a fibration if and only if it has the right lifting
property with respect to the maps in $A$.

\item A map is an acyclic fibration if and only if it has the
right lifting property with respect to the maps in $C$.

\item A map is an acyclic cofibration if and only if it has the
left lifting property with respect to the fibrations.
\end{enumerate}
Then there is a cofibrantly generated model category structure on
$\mathcal M$ in which $C$ is a set of generating cofibrations and
$A$ is a set of generating acyclic cofibrations.
\end{prop}

Let $\SSets$ denote the category of simplicial sets. Recall in
$\SSets$ we have for any $n \geq 0$ the $n$-simplex $\Delta [n]$,
its boundary $\dot \Delta [n]$, and, for any $0 \leq k \leq n$,
$V[n,k]$, which is $\dot \Delta [n]$ with the $k$th face removed.
Given a simplicial set $X$, we denote by $|X|$ its geometric
realization. The standard model category structure on $\SSets$ is
cofibrantly generated; the generating cofibrations are the maps
$\dot \Delta [n] \rightarrow \Delta [n]$  for $n \geq 0$, and the
generating acyclic cofibrations are the maps $V[n,k] \rightarrow
\Delta [n]$ for $n \geq 1$ and $0 \leq k \leq n$. More details on
simplicial sets and the model category structure on them can be
found in \cite{gj}.

There is a functor
\begin{equation}\label{U}
U: \SSets \rightarrow \mathcal{SCats}
\end{equation}
which takes a simplicial set $X$ to the category with objects $x$
and $y$ and with $\Hom (x,y)=X$ but no other nonidentity
morphisms.

We will say that a simplicial set $K$ is \emph{weakly
contractible} if all the homotopy groups of $|K|$ are trivial.

We will refer to the model category structure on the category of
simplicial categories with a fixed set $\mathcal O$ of objects,
denoted $\sco$, such that all the morphisms induce the identity
map on the objects, as defined by Dwyer and Kan in \cite{dk}.  The
weak equivalences are the maps satisfying condition W1 and the
fibrations are the maps satisfying condition F1.

We can then define our generating cofibrations and acyclic
cofibrations as follows.  The generating cofibrations are the maps
\begin{itemize}
\item (C1) $U \dot \Delta [n] \rightarrow U \Delta [n]$ for $n
\geq 0$, and

\item (C2) $\phi \rightarrow \{x\}$, where $\phi$ is the
simplicial category with no objects and $\{x\}$ denotes the
simplicial category with one object $x$ and no nonidentity
morphisms.
\end{itemize}
The generating acyclic cofibrations are
\begin{itemize}
\item (A1) the maps $U V[n,k] \rightarrow U \Delta [n]$ for $n
\geq 1$, and

\item (A2) inclusion maps $\{x\} \rightarrow \mathcal H$ which are
DK-equivalences, where $\{x\}$ is as in C2 and $\{\mathcal H\}$ is
a set of representatives for the isomorphism classes of simplicial
categories with two objects $x$ and $y$, weakly contractible
function complexes, and only countably many simplices in each
function complex.  Furthermore, we require that the inclusion map
$\{x\} \amalg \{y\} \rightarrow \mathcal H$ be a cofibration in
$\scxy$.
\end{itemize}

Notice that it follows from the definition of DK-equivalence that
the map $\{x\} \rightarrow \pi_0 \mathcal H$ is an equivalence of
categories.  In particular, all 0-simplices of the function
complexes of $\mathcal H$ are homotopy equivalences.

The idea behind the set A2 of generating acyclic cofibrations is
the fact that two simplicial categories can have a weak
equivalence between them which is not a bijection on objects, much
as two categories can be equivalent even if they do not have the
same objects.  We only require that our weak equivalences are
surjective on equivalence classes of objects.  Thus, we must
consider acyclic cofibrations for which the object sets are not
isomorphic.  The requirement that there be only countably many
simplices is included so that we have a set rather than a proper
class of maps.

We also show that this model structure is right proper. Recall
that a model category is \emph{right proper} if every pullback of
a weak equivalence along a fibration is a weak equivalence.  It
would be helpful to know if this model structure is also
\emph{left proper}, namely, that every pushout of a weak
equivalence along a cofibration is a weak equivalence.  We do not
know if this condition holds for $\mathcal{SC}$.  It would also be
useful to know that $\mathcal {SC}$ had the additional structure
of a simplicial model category, but we currently do not know of
such a structure.

In section 2, we show that these proposed generating acyclic
cofibrations satisfy the necessary conditions to be a generating
set. In section 3, we complete the proof of Theorem 1 and show
that the model structure is right proper. In section 4, we prove a
technical lemma that we needed in section 2.

\begin{thank}
The author would like to thank Bill Dwyer for helpful
conversations about this paper.
\end{thank}

\section{The Generating Acyclic Cofibrations}

In this section, we will show that a map in $\mathcal {SC}$
satisfies properties F1 and F2 if and only if it has the right
lifting property with respect to the maps in A1 and A2.

We begin by stating some facts about the Dwyer-Kan model category
structure $\sco$ on the category of simplicial categories with a
fixed set $\mathcal O$ of objects.  The weak equivalences are the
maps which satisfy property W1, and the fibrations are the maps
which satisfy property F1.  The cofibrations are the maps which
have the left lifting property with respect to the acyclic
fibrations.  However, we would like a more explicit description of
the cofibrations in this category, for which we need some
definitions.  If $\mathcal C$ is a simplicial category, then let
$\mathcal C_k$ denote the (discrete) category whose morphisms are
the $k$-simplices of the morphisms of $\mathcal C$.

\begin{definition}\cite[7.4]{dk}
A map $f:\mathcal C \rightarrow \mathcal D$ in $\sco$ is
\emph{free} if
\begin{enumerate}
\item $f$ is a monomorphism,

\item if $\ast$ denotes the free product, then in each simplicial
dimension $k$, the category $\mathcal D_k$ admits a unique free
factorization $\mathcal D_k= f(\mathcal C_k) \ast \mathcal F_k$,
where $\mathcal F_k$ is a free category, and

\item for each $k \geq 0$, all degeneracies of generators of
$\mathcal F_k$ are generators of $\mathcal F_{k+1}$.
\end{enumerate}
\end{definition}

\begin{definition}\cite[7.5]{dk}
A map $f:\mathcal C \rightarrow \mathcal D$ of simplicial
categories is a \emph{strong retract} of a map $f':\mathcal C
\rightarrow \mathcal D'$ if there exists a commutative diagram
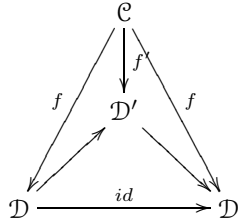
\[ \xymatrix{& \mathcal C \ar[ldd]_f \ar[rdd]^f \ar[d]^{f'}& \\
& \mathcal D' \ar[dr] & \\
\mathcal D \ar[ur] \ar[rr]^{id} && \mathcal D} \]
\end{definition}

Then, we have \cite[7.6]{dk} that the cofibrations of $\sco$ are
precisely the strong retracts of free maps.  In particular, a
cofibrant simplicial category is a retract of a free category.

Given these facts, we now continue with our discussion of the
generating acyclic cofibrations.

Recall~(\ref{U}) the map $U: \SSets \rightarrow \mathcal{SC}$. We
first consider the set A1 of maps $UV[n,k] \rightarrow U\Delta
[n]$ for $n \geq 1$ and $0 \leq k \leq n$. Using the model
category structure on simplicial sets, we can see that a map of
simplicial categories has the right lifting property with respect
to the maps in A1 if and only if it satisfies the property F1.

We then consider the maps in A2 which we would also like to be
generating acyclic cofibrations and show that maps with the right
lifting property with respect to the maps in A1 and A2 are
precisely the maps which satisfy conditions F1 and F2.  The proof
of this statement will take up the remainder of this section, and
we will treat each implication separately.

\begin{prop} \label{AAF}
Suppose that a map $f:\mathcal C \rightarrow \mathcal D$ of
simplicial categories has the right lifting property with respect
to the maps in A1 and A2.  Then $f$ satisfies condition F2.
\end{prop}

Before proving this proposition, we state a lemma.

\begin{lemma} \label{Factor}
Let $\mathcal F$ be a (discrete) simplicial category with object
set $\{x,y\}$ and one nonidentity morphism $g:x \rightarrow y$.
Let $\mathcal E'$ be a simplicial category also with object set
$\{x,y\}$. Let $i:\mathcal F \rightarrow \mathcal E'$ send $g$ to
a homotopy equivalence in $\Hom_{\mathcal E'}(x,y)$. This map $i$
can be factored as a composite $\mathcal F \rightarrow \mathcal H
\rightarrow \mathcal E'$ in such a way that the composite map
$\{x\} \rightarrow \mathcal F \rightarrow \mathcal H$ is
isomorphic to a map in A2.
\end{lemma}

We will prove this lemma in section 4.  We now prove Proposition
\ref{AAF} assuming Lemma \ref{Factor}.

\begin{proof}[Proof of Proposition \ref{AAF}]
Given objects $a_1$ in $\mathcal C$ and $b$ in $\mathcal D$, we
need to show that a homotopy equivalence $e:fa_1 \rightarrow b$ in
$\mathcal D$ lifts to a homotopy equivalence $d:a_1 \rightarrow
a_2$ for some $a_2$ in $\mathcal C$ such that $fa_2=b$ and $fd=e$.
So, we begin by considering the objects $a=fa_1$ and $b$ in
$\mathcal D$.

We first consider the case where $a \neq b$. Define $\mathcal E'$
to be the full simplicial subcategory of $\mathcal D$ with objects
$a=fa_1$ and $b$, and let $\mathcal F$ be a simplicial category
with objects $a$ and $b$ and a single nonidentity morphism $g:a
\rightarrow b$. Let $i:\mathcal F \rightarrow \mathcal E'$ send
$g$ to a homotopy equivalence $e:a \rightarrow b$.  By Lemma
\ref{Factor}, we can factor this map as $\mathcal F \rightarrow
\mathcal H \rightarrow \mathcal E'$ in such a way that the
composite $\{a\} \rightarrow \mathcal F \rightarrow \mathcal H$ is
isomorphic to a map in A2.

It follows that the composite $\{a_1\} \rightarrow \{a\}
\rightarrow \mathcal H$ is also isomorphic to a map in A2.  Then
consider the composite $\mathcal H \rightarrow \mathcal E'
\rightarrow \mathcal D$ where the map $\mathcal E' \rightarrow
\mathcal D$ is the inclusion map. These maps fit into a diagram
\[ \xymatrix{\{a_1\} \ar[r] \ar[d] & \mathcal C \ar[d]^f \\
\mathcal H \ar[r] \ar@{-->}[ur] & \mathcal D} \]

The lift exists because we assume that the map $f:\mathcal C
\rightarrow \mathcal D$ has the right lifting property with
respect to all maps in A2. Now, composing the map $\mathcal F
\rightarrow \mathcal H$ with the lift sends the map $g$ in
$\mathcal F$ to a map $d$ in $\mathcal C$ such that $fd=e$.  The
map $d$ is a homotopy equivalence since all the morphisms of
$\mathcal H$ are homotopy equivalences and therefore map to
homotopy equivalences in $\mathcal C$.

Now suppose that $a=b$.  Define $\mathcal E'$ to be the simplicial
category with two objects $a$ and $a'$ such that each function
complex of $\mathcal E'$ is the simplicial set $\Hom_\mathcal
D(a,a)$ and compositions are defined as they are in $\mathcal D$.
We then define the map $\mathcal E' \rightarrow \mathcal D$ which
sends both objects of $\mathcal E'$ to $a$ in $\mathcal D$ and is
the identity map on all the function complexes.  Given this
simplicial category $\mathcal E'$, the argument proceeds as above.
\end{proof}

We now prove the converse.

\begin{prop} \label{FFA}
Suppose $f:\mathcal C \rightarrow \mathcal D$ is a map of
simplicial categories which satisfies properties F1 and F2.  Then
$f$ has the right lifting property with respect to the maps in A2.
\end{prop}

Again, we state a lemma before proceeding with the proof of this
proposition.

\begin{lemma} \label{Lift}
Suppose that $A \rightarrow B$ is a cofibration, $C \rightarrow D$
is a fibration, and $B' \rightarrow B$ is a weak equivalence in a
model category $\mathcal M$.  Then in the following commutative
diagram
\[ \xymatrix{A \ar[r]^= \ar[d] & A  \ar[r] \ar@{^{(}->}[d] & C \ar@{->>}[d] \\
B' \ar[r]^\sim & B \ar[r] & D} \] a lift $B \rightarrow C$ exists
if and only if a lift $B' \rightarrow C$ exists.
\end{lemma}

\begin{proof}
If the lift $B \rightarrow C$ exists, it follows that the lift $B'
\rightarrow C$ exists via composition with the map $B' \rightarrow
B$.

To prove the converse, we first note that the map $B' \rightarrow
B$ can be factored as the composite
\[ B' \hookrightarrow B'' \twoheadrightarrow B \]
of a cofibration and a fibration, where each is a weak equivalence
because the map $B' \rightarrow B$ is.  Thus the fact that $A
\rightarrow B$ is a cofibration implies that there is a lift in
the diagram
\[ \xymatrix{A \ar[r] \ar@{^{(}->}[d] & B' \ar[r] & B'' \ar@{->>}[d]^\sim \\
B \ar@{-->}[urr] \ar[rr]^= && B} \]

It now suffices to show that there is a lift in the diagram
\[ \xymatrix{B' \ar[rr] \ar@{^{(}->}_\sim[d] && C \ar@{->>}[d] \\
B'' \ar[r] \ar@{-->}[urr] & B \ar[r] & D} \] However, this fact
again follows from the lifting properties of a model category.
\end{proof}

We are now able to prove the proposition.

\begin{proof}[Proof of Proposition \ref{FFA}]
We need to show that there exists a lift in any diagram of the
form
\[ \xymatrix{ \{x\} \ar[r] \ar[d] & \mathcal C \ar[d]^f \\
\mathcal H \ar[r] \ar@{-->}[ur] & \mathcal D} \] where $\{x\}
\rightarrow \mathcal H$ is a map in A2. Given an object $a_1$ in
$\mathcal C$ and a homotopy equivalence $e:a=fa_1 \rightarrow b$,
there exists an object $a_2$ in $\mathcal C$ and a homotopy
equivalence $d:a_1 \rightarrow a_2$ such that $fd=e$, since the
map $f$ satisfies property F2.

Let $g:x \rightarrow y$ be a homotopy equivalence in $\mathcal H$.
Let $\mathcal F$ denote the subcategory of $\mathcal H$ consisting
of the objects $x$ and $y$ and $g$ its only nonidentity morphism.
Consider the composite map $\{x\} \rightarrow \mathcal F
\rightarrow \mathcal H$ and the resulting diagram
\[ \xymatrix{\{x\} \ar[r] \ar[d] & \mathcal C \ar[dd] \\
\mathcal F \ar[d] \ar@{-->}[ur] & \\
\mathcal H \ar[r] & \mathcal D} \] Because the map $\mathcal F
\rightarrow \mathcal D$ factors through $\mathcal H$, which
consists of homotopy equivalences, the image of $g$ in $\mathcal
D$ is a homotopy equivalence.  Thus, the existence of the lift in
the above diagram follows from the fact that the map $f$ satisfies
F2.

Now, we need to show that the rest of $\mathcal H$ lifts to
$\mathcal C$.  We begin by assuming that $a \neq b$ and therefore
$a_1 \neq a_2$.  Consider the full simplicial subcategory of
$\mathcal C$ with objects $a_1$ and $a_2$, and denote by $\mathcal
C'$ the isomorphic simplicial category with objects $x$ and $y$.
Define $\mathcal D'$ analogously where we take objects $x$ and $y$
rather than $a$ and $b$.  Now, we can work in the category $\scxy$
of simplicial categories with fixed object set $\{x,y\}$. Note
that the map $\mathcal C' \rightarrow \mathcal D'$ is still a
fibration in $\mathcal {SC}$, and in fact it is a fibration in
$\scxy$. Now define $\mathcal E$ to be the pullback in the diagram
\[ \xymatrix{\mathcal E \ar[r] \ar[d] & \mathcal C' \ar[d] \\
\mathcal H \ar[r] & \mathcal D'} \] Then the map $\mathcal E
\rightarrow \mathcal H$ is also a fibration in $\scxy$
\cite[3.14(iii)]{ds}.

By Lemma \ref{Factor}, we can factor the map $\mathcal F
\rightarrow \mathcal E$ as the composite $\mathcal F \rightarrow
\mathcal H' \rightarrow \mathcal E$ for some simplicial category
$\mathcal H'$ such that the composite $\{x\} \rightarrow \mathcal
F \rightarrow \mathcal H'$ is isomorphic to a map in A2. Then,
note that the composite map $\mathcal H' \rightarrow \mathcal E
\rightarrow \mathcal H$ is a weak equivalence in $\scxy$ since all
the function complexes of $\mathcal H$ and $\mathcal H'$ are
weakly contractible.

Now, we have a diagram
\[ \xymatrix{\mathcal F \ar[r]^= \ar[d] & \mathcal F \ar[r] \ar[d] & \mathcal E \ar[d] \\
\mathcal H' \ar[r]^\sim & \mathcal H \ar[r] \ar@{-->}[ur] &
\mathcal D} \] in which the dotted arrow lift exists by Lemma
\ref{Lift}.

If $a=b$, then $\mathcal D'$ (and possibly $\mathcal C'$) as
defined above will have only one object $x$.  If this is the case,
then define the simplicial category $\mathcal D''$ with two
objects $x$ and $y$ such that each function complex is the
simplicial set $\Hom_{\mathcal D'}(x,x)$ (as in the proof of
Proposition \ref{AAF}).  We can then factor the map $\mathcal H
\rightarrow \mathcal D'$ through the object $\mathcal D''$, where
the map $\mathcal D'' \rightarrow \mathcal D$ sends both objects
of $\mathcal D''$ to $a$ in $\mathcal D$ and is the identity map
on each function complex. If $\mathcal C'$ also has one object,
then we obtain a simplicial category $\mathcal C''$ in the same
way. Then, we can repeat the argument above in the left-hand
square of the diagram
\[ \xymatrix{\mathcal F \ar[r] \ar[d] & \mathcal C'' \ar[r] \ar[d] &
\mathcal C' \ar[d] \\
\mathcal H \ar[r] \ar@{-->}[ur] & \mathcal D'' \ar[r] & \mathcal
D'} \] to obtain a lift $\mathcal H \rightarrow \mathcal C''$, and
hence a lift $\mathcal H \rightarrow \mathcal C'$ via composition.
\end{proof}

\section{The Model Category Structure}

In order to show that our proposed model category structure
exists, we need to show that our definitions are compatible with
one another. In particular, we need to prove that the maps with
the left lifting property with respect to the fibrations are
exactly the acyclic cofibrations, and that the maps with the right
lifting property with respect to the generating cofibrations are
exactly the maps which are fibrations and weak equivalences.
Before proving these statements, however, we prove that first
three model category axioms hold in $\mathcal{SC}$.

\begin{prop}
The category $\mathcal {SC}$ has all finite limits and colimits,
and its class of weak equivalences is closed under retracts and
satisfies the ``two out of three" property.
\end{prop}

\begin{proof}
It can be shown that the category of all simplicial categories has
all coproducts and all coequalizers, and therefore all finite
colimits, and all products and equalizers, and therefore all
finite limits. To prove the existence of coequalizers, for
example, we use the existence of coequalizers for sets (for the
objects) and simplicial sets (for the morphisms).  The two
properties for the class of weak equivalences follow as usual, for
example, as in \cite[8.10]{ds}.
\end{proof}

We first consider the sets C1 and C2. Suppose we have a map
$f:\mathcal C \rightarrow \mathcal D$ which is a fibration and a
weak equivalence.  Using the model category of simplicial sets, we
can see that a map satisfies conditions F1 and W1 if and only if
it has the right lifting property with respect to the maps $U \dot
\Delta [n] \rightarrow U \Delta [n]$ for $n \geq 0$, where $U$ is
the map~(\ref{U}) from simplicial sets to simplicial categories
defined in the first section.

However, the maps $U \dot \Delta [n] \rightarrow U \Delta [n]$
only generate those cofibrations between simplicial categories
with the same number of objects, a condition that we do not
require on our cofibrations of simplicial categories. Therefore,
we include as a generating cofibration the map $\phi \rightarrow
\{x\}$ from the simplicial category with no objects to the
single-object simplicial category with no nonidentity morphisms.
In other words, we are including the addition of an object as a
cofibration.

\begin{prop} \label{AFibs}
A map in $\mathcal{SC}$ is a fibration and a weak equivalence if
and only if it has the right lifting property with respect to the
maps in C1 and C2.
\end{prop}

\begin{proof}
First suppose that $f: \mathcal C \rightarrow \mathcal D$ is both
a fibration and a weak equivalence.  By conditions F1 and W1, the
map $\Hom_\mathcal C(a,b) \rightarrow \Hom_\mathcal D(fa,fb)$ is
an acyclic fibration of simplicial sets for any choice of objects
$a$ and $b$ in $\mathcal C$.  In other words, there is a lift in
any diagram of the form
\[ \xymatrix{\dot \Delta [n] \ar[r] \ar[d] & \Hom_\mathcal C(a,b)
\ar[d] \\
\Delta [n] \ar[r] \ar@{-->}[ur] & \Hom_\mathcal D(fa,fb)} \]
However, having this lift is equivalent to having a lift in the
diagram
\[ \xymatrix{U \dot \Delta [n] \ar[r] \ar[d] & \mathcal C \ar[d] \\
U \Delta [n] \ar[r] \ar@{-->}[ur] & \mathcal D} \] where the
objects $x$ and $y$ of $U \dot \Delta [n]$ map to $a$ and $b$ in
$\mathcal C$, and analogously for $U \Delta [n]$ and $\mathcal D$.
Hence, $f$ has the right lifting property with respect to the maps
in C1.

It remains only to show that $f$ has the right lifting property
with respect to the map $\phi \rightarrow \{x\}$. However, this
property is equivalent to $f$ being onto on objects. Being onto on
homotopy equivalence classes of objects follows from condition W2.
Then suppose that $e:a \rightarrow b$ is an isomorphism in
$\mathcal D$ and there is an object $a_1$ in $\mathcal C$ such
that $fa_1 =a$. Since $e$ is a homotopy equivalence, by F2 there
is a homotopy equivalence in $\mathcal C$ with domain $a_1$ and
which maps to $e$ under $f$. In particular, there is an object
$a_2$ in $\mathcal C$ mapping to $b$.

Conversely, suppose that $f$ has the right lifting property with
respect to the maps in C1 and C2.  Again, using the model category
structure on simplicial sets, we have that the map
\[ \Hom_\mathcal C (a,b) \rightarrow \Hom_\mathcal D (fa,fb) \]
is both a fibration and a weak equivalence, satisfying both F1 and
W1. It follows that $\Hom_{\pi_0 \mathcal C}(a,b) \rightarrow
\Hom_{\pi_0 \mathcal D}(fa,fb)$ is an isomorphism.  As above,
having the right lifting property with respect to the map $\phi
\rightarrow \{x\}$ is equivalent to being onto on objects. These
two facts show then that $\pi_0 \mathcal C \rightarrow \pi_0
\mathcal D$ is an equivalence of categories, proving condition W2.

It remains to show that $f$ satisfies property F2. By Proposition
\ref{AAF} and the fact that satisfying F1 is equivalent to having
the right lifting property with respect to maps in A1, it suffices
to show that $f$ has the right lifting property with respect to
the maps in A2. But, a map $\{x\} \rightarrow \mathcal H$ in A2
can be written as a (possibly infinite) composition of a pushout
along $\phi \rightarrow \{x\}$ followed by pushouts along maps of
the form $U \dot \Delta [n] \rightarrow U \Delta [n]$, and $f$ has
the right lifting property with respect to all such maps since
these are just the maps in C1 and C2.
\end{proof}

\begin{prop} \label{ACofs}
A map in $\mathcal{SC}$ is an acyclic cofibration if and only if
it has the left lifting property with respect to the fibrations.
\end{prop}

The proof will require the use of the following lemma:

\begin{lemma} \label{Push}
Let $\mathcal A \rightarrow \mathcal B$ be a map in A1 or A2 and
$i:\mathcal A \rightarrow \mathcal C$ any map in $\mathcal {SC}$.
Then in the pushout diagram
\[ \xymatrix{\mathcal A \ar[r]^i \ar[d] & \mathcal C \ar[d] \\
\mathcal B \ar[r] & \mathcal D} \] the map $\mathcal C \rightarrow
\mathcal D$ is a weak equivalence.
\end{lemma}

\begin{proof}
First suppose that the map $\mathcal A \rightarrow \mathcal B$ is
in A2. Let $\mathcal O$ be the set of objects of $\mathcal C$ and
define $\mathcal O'$ to be the set $\mathcal O \backslash \{x\}$.
(For simplicity of notation, we assume that $ix=x$.)  Assume as
before that $x$ and $y$ are the objects of $\mathcal H$. We denote
also by $\mathcal O'$ the (simplicial) category with object set
$\mathcal O'$ and no nonidentity morphisms. Consider the diagram
\[ \xymatrix{\mathcal X =\mathcal O \amalg \{y\} \ar[r] \ar[d] & \mathcal C
\amalg \{y\}= \mathcal C' \ar[d] \\
\mathcal H' =\mathcal O' \amalg \mathcal H \ar[r] & \mathcal D} \]
and notice that $\mathcal D$ is also the pushout of this diagram.
Since $\mathcal X$ (regarded as a set) is the object set of any of
these categories, note that the left hand vertical arrow is a
cofibration in $\mathcal {SC}_\mathcal X$.

We factor the map $\mathcal X \rightarrow \mathcal C'$ as the
composite of a cofibration and an acyclic fibration in $\mathcal
{SC}_\mathcal X$
\[ \xymatrix@1{\mathcal X \ar@{^{(}->}[r] & \mathcal C''
\ar@{->>}[r]^\sim & \mathcal C'.} \]

Since $\mathcal {SC}_\mathcal X$ is proper \cite[7.3]{dk}, it
follows from \cite[13.5.4]{hirsch} that the pushouts of each row
in the diagram
\[ \xymatrix{\mathcal H' & \mathcal X \ar@{_{(}->}[l]
\ar[r] & \mathcal C' \\
\mathcal H' \ar[d]^\sim \ar[u]_= & \mathcal X \ar@{_{(}->}[l]
\ar@{^{(}->}[r] \ar[d]^= \ar[u]_= & \mathcal C'' \ar[d]^= \ar[u]_\sim \\
\pi_0 \mathcal H' & \mathcal X \ar[l] \ar@{^{(}->}[r] & \mathcal
C''} \] are weakly equivalent to one another.  In particular, the
pushout of the bottom row is weakly equivalent to $\mathcal D$. It
remains to show that there is a weak equivalence of pushouts of
the rows of the diagram
\[ \xymatrix{\pi_0 \mathcal H' \ar[d] & \mathcal X \ar[l] \ar[d] \ar@{^{(}->}[r] & \mathcal
C'' \ar[d] \\
\pi_0 \mathcal H' & \mathcal X \ar[r] \ar[l] & \mathcal C'.}
\] However, a calculation shows that the pushout of this bottom row is weakly equivalent in
$\mathcal {SC}$ to the pushout of the diagram
\[ \xymatrix@1{\pi_0 \mathcal H & \{x\} \ar[l] \ar[r] & \mathcal
C} \] and therefore that the pushout of the top row is weakly
equivalent to the pushout of the bottom row. It follows that the
map $\mathcal C \rightarrow \mathcal D$ is a weak equivalence in
$\mathcal {SC}$.

For the maps in A1, we have pushout diagrams
\[ \xymatrix{UV[n,k] \ar[r]^-j \ar[d] & \mathcal C \ar[d] \\
U\Delta [n] \ar[r] & \mathcal D}. \]  As before, define $\mathcal
O$ to be the object set of $\mathcal C$.  Let $\mathcal O'' =
\mathcal O \backslash \{x,y\}$.  (Again, for notational simplicity
we will assume that $jx=x$ and $jy=y$.)  Now we consider the
diagram
\[ \xymatrix{\mathcal O'' \amalg UV[n,k] \ar[r] \ar[d] & \mathcal
C \ar[d] \\
\mathcal O'' \amalg U\Delta [n] \ar[r] & \mathcal D} \] in $\sco$.
However, since the left vertical map is a weak equivalence and
assuming that the top map is a cofibration (factoring if necessary
as above), we can again use the fact that $\sco$ is proper to show
that $\mathcal C \rightarrow \mathcal D$ is a weak equivalence in
$\sco$ and thus also in $\mathcal {SC}$.
\end{proof}

\begin{proof}[Proof of Proposition \ref{ACofs}]
First suppose that a map $\mathcal C \rightarrow \mathcal D$ is an
acyclic cofibration.  By the small object argument (\cite[Sec.
7]{ds} or \cite[Ch. 11]{hirsch}), we have a factorization of the
map $\mathcal C \rightarrow \mathcal D$ as the composite $\mathcal
C \rightarrow \mathcal C' \rightarrow \mathcal D$ where $\mathcal
C'$ is obtained from $\mathcal C$ by a directed colimit of
iterated pushouts along the maps in A1 and A2.  Thus, by Lemma
\ref{Push} above and the fact that a directed colimit of such maps
is a weak equivalence, this map $\mathcal C \rightarrow \mathcal
C'$ is a weak equivalence.  Furthermore, the map $\mathcal C'
\rightarrow \mathcal D$ has the right lifting property with
respect to the maps in A1 and A2. Thus, by Proposition \ref{AAF},
it is a fibration.  It is also a weak equivalence since the maps
$\mathcal C \rightarrow \mathcal D$ and $\mathcal C \rightarrow
\mathcal C'$ are, by axiom MC2.  In particular, by the definition
of cofibration, it has the right lifting property with respect to
the cofibrations. Therefore, there exists a dotted arrow lift in
the diagram
\[ \xymatrix{\mathcal C \ar[r]^\sim \ar[d] & \mathcal C'
\ar[d] \\
\mathcal D \ar[r]^= \ar@{-->}[ur] & \mathcal D} \] Hence the map
$\mathcal C \rightarrow \mathcal D$ is a retract of the map
$\mathcal C \rightarrow \mathcal C'$ and therefore also has the
left lifting property with respect to fibrations.

Conversely, suppose that the map $\mathcal C \rightarrow \mathcal
D$ has the left lifting property with respect to fibrations.  In
particular, it has the left lifting property with respect to the
acyclic fibrations, so it is a cofibration by definition.  We
again obtain a factorization of this map as the composite
$\mathcal C \rightarrow \mathcal C' \rightarrow \mathcal D$ where
$\mathcal C'$ is obtained from $\mathcal C$ by iterated pushouts
of the maps in A1 and A2.  Once again, the map $\mathcal C'
\rightarrow \mathcal D$ has the right lifting property with
respect to the maps in A1 and A2 and thus is a fibration by
Proposition \ref{AAF}. Therefore there is a lift in the diagram
\[ \xymatrix{\mathcal C \ar[r] ^\sim \ar[d] & \mathcal C' \ar[d]
\\
\mathcal D \ar[r]^= \ar@{-->}[ur] & \mathcal D} \]  Again using
Lemma \ref{Push}, the map $\mathcal C \rightarrow \mathcal D$ is a
weak equivalence because it is a retract of the map $\mathcal C
\rightarrow \mathcal C'$.
\end{proof}

We have now proved everything we need for the existence of the
model category structure on $\mathcal {SC}$.

\begin{proof}[Proof of Theorem \ref{MC}]
It remains to show that the four conditions of Proposition
\ref{CofGen} are satisfied.  It can be shown that both $\phi$ and
$\{x\}$ are small, and using the smallness of $V[n,k]$ and $\Delta
[n]$ in $\SSets$ \cite[3.1.1]{hovey}, it can be shown that each $U
\dot \Delta [n]$ is small relative to the set C1 and each $U
V[n,k]$ is small relative to the set A1 \cite[10.5.12]{hirsch}.
Therefore, condition 1 holds.  Condition 2 follows from
Propositions \ref{AAF} and \ref{FFA}. Condition 3 is proved in
Proposition \ref{AFibs}, and condition 4 is proved in Proposition
\ref{ACofs}.
\end{proof}

We conclude this section with the following result about this
model category structure.

\begin{prop}
The model category structure $\mathcal{SC}$ is right proper.
\end{prop}

\begin{proof}
Suppose that
\[ \xymatrix{A=B \times_DC \ar[r]^-f \ar[d] & B \ar[d]^g \\
C \ar[r]^h & D} \] is a pullback diagram, where $g:B \rightarrow
D$ is a fibration and $h:C \rightarrow D$ is a DK-equivalence.  We
would like to show that $f:A \rightarrow B$ is a DK-equivalence.

We first need to show that $\Hom_A(x,y) \rightarrow \Hom_B(fx,fy)$
is a weak equivalence of simplicial sets for any objects $x$ and
$y$ of $A$.  However, this fact follows since the model category
structure on simplicial sets is right proper
\cite[13.1.4]{hirsch}.

It remains to prove that $\pi_0A \rightarrow \pi_0B$ is an
equivalence of categories.  After applying $\pi_0$ to the result
of the previous paragraph, it suffices to show that $A \rightarrow
B$ is essentially surjective on objects.

Consider an object $b$ of $B$ and its image $g(b)$ in $D$.  Since
$C \rightarrow D$ is a DK-equivalence, there exists some object
$c$ of $C$ together with a homotopy equivalence $g(b) \rightarrow
h(c)$ in $D$.  Since $B \rightarrow D$ is a fibration, there
exists an object $b'$ and homotopy equivalence $b \rightarrow b'$
in $B$ such that $g(b')=h(c)$.  Using the fact that $A$ is a
pullback, we have a homotopy equivalence $b \rightarrow
f((b',c))$, completing the proof.
\end{proof}

\section{Proof of Lemma \ref{Factor}}

Recall that we have a (simplicial) category $\mathcal F$ with
objects $x$ and $y$ and a single nonidentity morphism $g:x
\rightarrow y$, and a simplicial category $\mathcal E'$ also with
objects $x$ and $y$ such that there is a map $i:\mathcal F
\rightarrow \mathcal E'$ which sends $g$ to a homotopy equivalence
$x \rightarrow y$ in $\mathcal E'$. We first replace $\mathcal E'$
by its subcategory of homotopy equivalences which we denote by
$\mathcal E$. In order to make our constructions homotopy
invariant, we take functorial cofibrant replacements $\widetilde
{\mathcal F} \rightarrow \mathcal F$ and $\widetilde {\mathcal E}
\rightarrow \mathcal E$ in the model category $\scxy$ as given in
\cite[2.5]{dk}, and in this construction $\widetilde{\mathcal F}$
is actually isomorphic to $\mathcal F$.

Now, take the localization $\mathcal F^{-1} \mathcal F$
(respectively $\widetilde {\mathcal E}^{-1} \widetilde {\mathcal
E}$) obtained by formally inverting all the morphisms in each
simplicial degree of $\mathcal F$ (respectively $\widetilde
{\mathcal E}$). These localizations are the groupoid completions
of $\mathcal F$ and $\widetilde{\mathcal E}$, respectively. (In
taking a functorial cofibrant replacement and then the groupoid
completion, we have taken the simplicial localizations of
$\mathcal F$ and $\mathcal E$ with respect to all the morphisms in
each as defined in \cite{dk}.)  We now have a diagram
\[ \xymatrix{\mathcal F \ar[d] & \mathcal F \ar[r]
\ar[d] \ar[l]_= & \mathcal F^{-1} \mathcal F \ar[d] \\
\mathcal E & \widetilde{\mathcal E} \ar[l] \ar[r] &
\widetilde{\mathcal E}^{-1} \widetilde{\mathcal E}.} \]

To assure that our next step is homotopy invariant, we factor the
map $\mathcal F^{-1} \mathcal F \rightarrow \widetilde{\mathcal
E}^{-1} \widetilde{\mathcal E}$ as the composite
\[ \xymatrix@1{\mathcal F^{-1} \mathcal F
\ar[r]^-i & \mathcal Z \ar[r]^-p & \widetilde{\mathcal E}^{-1}
\widetilde{\mathcal E}} \] where $i$ is an acyclic cofibration and
$p$ is a fibration in $\scxy$.  However, to avoid more notation
than necessary, we will assume that the map $\mathcal F^{-1}
\mathcal F \rightarrow \widetilde{\mathcal E}^{-1}
\widetilde{\mathcal E}$ is a fibration and continue to write
$\mathcal F^{-1} \mathcal F$ rather than $\mathcal Z$.

We take the pullback of the bottom right hand corner of the above
diagram and denote it $\mathcal G$:
\[ \xymatrix{\mathcal G \ar[r] \ar[d] & \mathcal F^{-1} \mathcal F \ar[d] \\
\widetilde{\mathcal E} \ar[r] & \widetilde{\mathcal E}^{-1}
\widetilde{\mathcal E}.} \]  Notice that this diagram is a
pullback in $\scxy$ as well as in $\mathcal{SC}$.

\begin{lemma}
The composite map $\{x\} \rightarrow \mathcal F \rightarrow
\mathcal G$ is a weak equivalence in $\mathcal {SC}$.
\end{lemma}

\begin{proof}
Since the simplicial categories $\widetilde{\mathcal E}$,
$\widetilde{\mathcal E}^{-1} \widetilde{\mathcal E}$, and
$\mathcal F^{-1} \mathcal F$ all consist of homotopy equivalences,
so must $\mathcal G$. Therefore, all the morphisms of $\pi_0
\mathcal G$ are isomorphisms, and in particular, the objects $x$
and $y$ are isomorphic in $\pi_0 \mathcal G$.

It then suffices to show that $\mathcal G$ has weakly contractible
function complexes.  Because all of the morphisms of $\mathcal E$,
and hence also of $\widetilde{\mathcal E}$, are homotopy
equivalences, the map $\widetilde{\mathcal E} \rightarrow
\widetilde{\mathcal E}^{-1} \widetilde{\mathcal E}$ is a weak
equivalence in $\scxy$ \cite[9.5]{dk}.

Note that $\mathcal F^{-1} \mathcal F$ is the simplicial category
in $\scxy$ with exactly one morphism between any two objects. In
particular, $\mathcal F^{-1} \mathcal F$ has weakly contractible
function complexes.

Now, because all the categories have as objects $x$ and $y$ and
all the maps involved are the identity on these objects, we can
consider the above pullback diagram in $\scxy$.  Since this model
category structure is right proper \cite[7.3]{dk}, every pullback
of a weak equivalence along a fibration is a weak equivalence. The
map $\widetilde{\mathcal E} \rightarrow \widetilde{\mathcal
E}^{-1} \widetilde{\mathcal E}$ is a weak equivalence and the map
$\mathcal F^{-1} \mathcal F \rightarrow \widetilde{\mathcal
E}^{-1} \widetilde{\mathcal E}$ is a fibration, so it follows that
the map $\mathcal G \rightarrow \mathcal F^{-1} \mathcal F$ is a
weak equivalence in $\scxy$, and therefore $\mathcal G$ has weakly
contractible function complexes. Thus, the map $\{x\} \rightarrow
\mathcal G$ satisfies the conditions to be a weak equivalence in
$\mathcal {SC}$.
\end{proof}

However, not all the maps $\{x\} \rightarrow \mathcal G$ are
isomorphic to maps in A2 because the simplicial categories
$\mathcal G$ could have an uncountable number of simplices in
their function complexes. Furthermore, there is no reason to
assume that the inclusion map $\{x\} \amalg \{y\} \rightarrow
\mathcal G$ is a cofibration in $\scxy$. To complete the proof, we
need to show that any acyclic cofibration $\{x\} \rightarrow
\mathcal G$ as above factors as a composite $\{x\} \rightarrow
\mathcal H \rightarrow \mathcal G$ where the inclusion map $\{x\}
\rightarrow \mathcal H$ is in A2.

Let $\mathcal H_0$ be the simplicial category $\mathcal F$. Let
$i:\mathcal H_0 \rightarrow \mathcal G$ be the inclusion map. We
will construct a simplicial category $\mathcal H$ from $\mathcal
H_0$ satisfying the necessary properties specified in A2.  We
first state the following lemma:

\begin{lemma} \label{AB}
Let $f:A \rightarrow B$ be a map of simplicial sets where $B$ is
weakly contractible, and let $u:S^n \rightarrow |A|$ be a map of
CW-complexes for some $n \geq 0$.  Then $f$ can be factored as a
composite $A \rightarrow A' \rightarrow B$ where $A'$ is obtained
from $A$ by attaching a finite number of nondegenerate simplices
and the composite map of spaces $S^n \rightarrow |A| \rightarrow
|A'|$ is null homotopic.
\end{lemma}

\begin{proof}
We first assume that the map $f$ is a cofibration; if not, we
factor it as the composite
\[ \xymatrix@1{A \ar[r]^i & B' \ar[r]^p & B} \]
where in the model category structure on simplicial sets $i$ is a
cofibration and $p$ is an acyclic fibration. Thus, we can assume
that $f$ is an inclusion map, replacing $B$ by $B'$ if needed.

Now consider the composite map of spaces $S^n \rightarrow |A|
\rightarrow |B|$, which is necessarily null homotopic since $B$ is
weakly contractible.  The composite map then factors through
$CS^n$, the cone on $S^n$, and we have a diagram
\[ \xymatrix{S^n \ar[d] \ar[r] & |A| \ar[r] & |B| \\
CS^n \ar[urr] && } \]

Now, since $CS^n$ is compact, its image will intersect only a
finite number of cells of $|B|$ nontrivially.  Then define $A'$ to
be a simplicial set such that $|A'|$ contains $|A|$ as well as all
the cells in this image.
\end{proof}

Now, consider the categories $\mathcal H_0$ and $\mathcal G$ as
described above and the inclusion map $i:\mathcal H_0 \rightarrow
\mathcal G$. Each of these categories has four function complexes
to consider. For the category $\mathcal H_0$ we call them $H_j$,
and for $\mathcal G$ we call them $G_j$ for $1 \leq j \leq 4$.
(The numbering is arbitrary but must match up between the two
categories.  So if $H_1= \Hom_{\mathcal H} (x,y)$, then we must
have $G_1 = \Hom_\mathcal G(x,y)$.)

First notice that $\Hom_{\mathcal H_0}(y,x)=\phi$.  We begin by
adding a map $y \rightarrow x$ to make this function complex
nonempty, as well as all composites generated by it.  Then,
considering all four function complexes, we identify $n \geq 0$
such that all maps $S^m \rightarrow |H_j|$ are null homotopic for
all $0 \leq m < n$ and all $1 \leq j \leq 4$, but there is a map
$S^n \rightarrow |H_j|$ which is not null homotopic for some $j$.
We then apply Lemma \ref{AB} to the map $H_j \rightarrow G_j$ and
the map $S^n \rightarrow |H_j|$.

Replace the function complex $H_j$ with the simplicial set $A'$
obtained from Lemma \ref{AB}. This process may result in more maps
$S^m \rightarrow |A'|$ which are not null homotopic than for the
original $|H_j|$, but only for $m>n$. Also, it will not have more
than countably many more such maps than $|H_j|$ did. Now that we
have added simplices to our function complex, we include all
necessary compositions of these morphisms with the original
morphisms of $\mathcal H_0$ to obtain a new simplicial category
which we denote $\mathcal H_1$. There will be at most countably
many new simplices added from these compositions. Repeat the above
process with another map from $S^n$ to a function complex of
$\mathcal H_1$, again, where $n$ is minimal, to obtain another
category $\mathcal H_2$. Continue, perhaps countably many times,
to obtain a category $\mathcal H$ such that for any $n$ and any
function complex $H'$ of $\mathcal H$, any map $S^n \rightarrow
|H'|$ is nullhomotopic. To show that it is possible to obtain such
an $\mathcal H$ in this way, we need only show that there are at
most countably many homotopy classes of maps from spheres to each
function complex that need to be killed off. However, this fact
follows from the following lemma:

\begin{lemma} \label{Count}
Let $A$ be a simplicial set with countably many simplices. Then
for all $n \geq 0$ there are at most countably many distinct
homotopy classes of maps $S^n \rightarrow |A|$.
\end{lemma}

\begin{proof}
It suffices to show that there are at most countably many homotopy
classes of maps from $S^n$ into any finite CW complex $X$.  For a
simply connected CW complex $X$, an argument using Serre mod
$\mathcal C$ theory \cite{serre} shows that all the homotopy
groups of $X$ are countable if and only if the homology groups of
$X$ are countable, which they are when $X$ is finite.  The case of
a general CW complex $X$ follows from this one using a universal
cover argument.
\end{proof}

By construction, this simplicial category $\mathcal H$ is free,
and therefore the map
\[ \{x\} \amalg \{y\} \rightarrow \mathcal H \] is a cofibration in
$\scxy$. Thus, we have obtained a factorization $\{x\} \rightarrow
\mathcal H \rightarrow \mathcal G$.  We are now able to complete
the proof of Lemma \ref{Factor}.

\begin{proof}[Proof of Lemma \ref{Factor}]
Using the simplicial category $\mathcal H$ from above and the map
$\mathcal H \rightarrow \mathcal G$, we obtain a composite map
\[ \{x\} \rightarrow \mathcal F \rightarrow \mathcal H \rightarrow \mathcal G
\rightarrow \widetilde{\mathcal E} \rightarrow \mathcal E
\rightarrow \mathcal E'.\] In particular, we have a factorization
$\mathcal F \rightarrow \mathcal H \rightarrow \mathcal E'$.  As
we have shown above, the composite $\{x\} \rightarrow \mathcal F
\rightarrow \mathcal H$ is isomorphic to a map in A2.
\end{proof}

\end{document}